\documentclass{article}
\usepackage{amsmath}
\usepackage{amssymb}
\usepackage{amscd}

\def\be{\begin{equation}}       \def\ee{\end{equation}}
\def\bd{\begin{displaymath}}    \def\ed{\end{displaymath}}
\def\beq{\begin{eqnarray}}      \def\eeq{\end{eqnarray}}
\def\bseq{\begin{eqnarray*}}    \def\eseq{\end{eqnarray*}}
\def\ba{\begin{array}}          \def\ea{\end{array}}
\def\ben{\begin{enumerate}}     \def\een{\end{enumerate}}
\def\nonum{\nonumber}

\newtheorem{Prop}{Proposition}
\def\bpr{\begin{Prop}}           \def\epr{\end{Prop}}

\def\ra{\rightarrow}
\def\cp{\circ}
\def\lam{\Omega}
\def\lop{\hat{\Omega}\,}
\def\qop{\hat{Q}\,}
\def\hop{\hat{H}\,}
\def\deropx{~_{\Omega}D_x}
\def\derop{~ _{\Omega}D}
\def\Elam{E_{\Omega}}
\def\derinv{~ _{\Omega}D^{-1}_x}
\def\lopr{\hat{\Omega}_{{\lambda},k}\,}

\def\xop{\hat{x}\,}
\def\dx{{\partial}_x}
\def\Ellam{E_{\lambda}}
\def\lopf{\hat{\Omega}_{{\lambda},\mu}\,}
\def\Lmn{\hat{\Omega}_{\lambda , \mu}\,}
\def\loprinv{\hat{\Omega}^{-1}_{{\lambda},k}\,}
\def\derlr{~_{{\Omega}_{\lambda, k}}D_x}

\numberwithin{equation}{section}
\begin{document}
\begin{center}
{\large\bf CONSTRUCTION OF THE GENERALISED $q$--DERIVATIVE OPERATORS}\\
\bigskip \bigskip \bigskip 
\bigskip
\bigskip
\bigskip
\bigskip

{\large\bf Dayanand Parashar}\\
{\sl Department of Physics, A.R.S.D. College, University of Delhi\\ 
New Delhi - 110021, India}\\
\bigskip \bigskip
{\large\bf Deepak Parashar}\\
{\sl Department of Mathematics, University of Wales
Swansea, Singleton Park\\ 
Swansea SA2 8PP, United Kingdom}\\
\smallskip
{\normalsize\tt D.Parashar@swansea.ac.uk\\
http://www-maths.swan.ac.uk/staff/dp/}\\
\bigskip \bigskip \bigskip \bigskip
\bigskip \bigskip
\bigskip
\bigskip
{\large \bf Abstract}\\
\end{center}
\medskip

This investigation pertains to the construction of a class of generalised
deformed derivative operators which furnish the familiar finite difference
and the $q$--derivatives as special cases. The procedure involves the
introduction of a linear operator which is multiplicative over functions
of a real variable. The validity of the general prescriptions is
ascertained by considering suitable examples of such derivatives and
constructing their eigenfunctions explicitly. The relationship of a
particular version of the operator with the one-dimensional M\"obius
transformation is also established.\\

MSC: {\tt 05A30; 26A24; 81R50}
\bigskip \bigskip \bigskip \bigskip
\bigskip \bigskip \bigskip \bigskip
\begin{center}
{\sl J. Geom. Phys. 48 (2003) 297--308}
\end{center}

\newpage

\section{Introduction}

Ever since the pioneering work of Jackson \cite{jack}, there have been
numerous attempts \cite{coon, ext, agar} to carry out investigation into
various aspects of deformations which have played a crucial role in the
understanding of important mathematical and physical concepts. A concrete
realisation of such attributes is easily traced to the emergence of
quantum groups \cite{drin, frt}, embodying appropriate deformations of
algebraic and geometric structures in a particularly elegant
formulation. While there are several ways in which to study deformations,
perhaps the most essential ingredients in such a study would invariably
comprise, for instance, the $q$--deformed derivative operators, coordinate
multiplication operators, etc., among others. Significantly, these
$q$--deformed operators have met with tremendous success when applied to a
host of problems, most notably in the realm of oscillator algebra
\cite{coon, arik, bied} and $q$--calculus \cite{ks}. This offers the main
motivation to seek a generalisation of these operators in order to be able
to try applications of the underlying prescriptions to a wider class of
systems.
\par
In Sec. $2$, we construct a generalised deformed derivative operator and
study the corresponding eigenfunctions. Sec. $3$ is devoted to a
discussion of specific examples of the generalised version of this
operator where we explicitly demonstrate the salient features for the
purpose of applications. A possible connection of the deformed operator
with the reduced M\"obius transformation operator is discussed in
Sec. $4$. Concluding remarks form contents of Sec. $5$.

\section{The Deformed Derivative}

To facilitate construction of the generalised deformed derivative, we
consider a linear operator $\lop$ whose action on the product of two
functions, say $f(x)$ and $g(x)$, $x\in \mathbb{R}$, can be written as
\begin{equation}
 \lop \{ f(x)g(x) \} = \{ \lop f(x) \}\{\lop g(x) \}
  = f(\lop x )g( \lop x )
\end{equation}
clearly displaying the property that $\lop$ is multiplicative over
functions of $x$. We further assign to $\lop$ dependence on a parameter
$\lambda$ such that $ \lop \rightarrow I$ (the identity operator) in the
limit $\lambda \rightarrow 0$. This requirement is necessitated by the
desire to accomplish meaningful comparisons with the more familiar
derivative operators. In particular, the translation operator
$\Omega = \hop = e^{h \partial _x}$ such that $\hop f(x) = f(x+h)$, while
the $q$--dilation operator  $\Omega = \qop = e^{\ln(q) x \dx}$ such that
$\qop f(x) = f(q x)$ form examples of such an operator. An appropriate 
generalisation of the deformed derivative operator is constructed through
the definition
\begin{equation}\label{derop}
  \deropx =
  \left\{ \frac{1}{(\lop -I)x} \right\} \left( \lop - I \right)
\end{equation}
which, operating on any function $f(x)$, gives
\begin{equation}\label{deropf}
  \deropx f(x) = \frac{ f(\lop x) - f(x)}{(\lop x) - x}
\end{equation}
It follows immediately that $\deropx x = 1$ and $\deropx c = 0$
for any constant $c$. 
\begin{Prop}
The deformed derivative $\deropx$ obeys two modified
Leibniz rules
\begin{equation}\label{Leib1}
  \deropx \{f(x) g(x) \} =
  \{ \deropx f(x) \} g(x) + \{ \lop f(x) \} \deropx g(x)
\end{equation}
\begin{equation}\label{Leib2}
  \deropx \{f(x) g(x) \} =
  f(x) \{ \deropx g(x) \} + \{ \deropx f(x) \} \lop g(x)
\end{equation}
and the chain rule
\begin{equation}
\deropx f(g(x)) = {\derop}_{g(x)} f(g(x)) \deropx g(x)
\end{equation}
\end{Prop}
{\bf Proof:}
\begin{eqnarray*}
\deropx \{f(x) g(x) \}
&=& \frac{1}{(\lop -I)x}\left\{ \lop f(x)g(x) - f(x)g(x)\right\} \\
&=& \frac{1}{(\lop -I)x} \{ [f(\lop x)]g(x) - [f(\lop x)]g(x) \\
& &  \hfill   +  f(\lop x)g(\lop x) - f(x)g(x)     \} \\
&=& \frac{1}{(\lop -I)x} \left\{ [f(\lop x) - f(x)]g(x) +
                   f(\lop x)[g(\lop x) - g(x)]     \right\} \\
&=& \left\{\frac{ f(\lop x) - f(x)}{(\lop x) - x}\right\}g(x)
 + f(\lop x)\left\{\frac{ g(\lop x) - g(x)}{(\lop x) -  x}\right\}\\
&=& \left\{ \deropx f(x) \right\} g(x)
        + \lop f(x) \left\{ \deropx g(x) \right\}
\end{eqnarray*}
The second Leibniz rule (\ref{Leib2}) follows from the first if we
simply interchange the roles of $f$ and $g$. To prove the modified chain
rule, we have
\begin{eqnarray*}
\deropx \{f\circ g\} (x)
&=& \frac{ \{f\circ g\} (\lop x) - \{f\circ g\}(x)}{(\lop x) - x}\\
&=& \left\{
        \frac{f(g(\lop x)) - f(g(x))}{g(\lop x)- g(x)}
    \right\}
    \left\{
    \frac{g(\lop x) - g(x)}{(\lop x)- x}
    \right\} \\
&=& \left\{
        \frac{(\lop - I)f(g(x))}{(\lop - I)g(x)}
    \right\}
    \left\{
    \frac{(\lop - I)g(x)}{(\lop - I)x}
    \right\} \\
&=& {\derop}_{g(x)} f(g(x)) \deropx g(x)
\end{eqnarray*}\\
In the limit $\lambda \rightarrow 0$, we see that $\deropx \rightarrow
\frac{d}{dx}=\dx$. Note that this generalised derivative specialises to
the familiar finite difference derivative $~_{h}\Delta_x$:
  \[ ~_{h}\Delta_x f(x) = \frac{f(x+h)-f(x)}{h} \]
and the $q$--derivative $~_{q}D_x$:
  \[ ~_{q}D_x f(x) = \frac{f(qx)-f(x)}{(q-1)x} \]
as expected.

\subsection{Inverse deformed derivative ($\derinv$)}

The nonsingular behaviour of the operator $\deropx$ can be ascertained
only if the existence of its inverse $\derinv$ is guaranteed. Consequently, 
the condition $\deropx \cp \derinv = I$ must necessarily be satisfied on
analytic functions. We, therefore, proceed to define such an inverse. The
operation of the deformed derivative $\deropx$ on an arbitrary analytic
function $f(x)$ can be rewritten in the form
\begin{equation}
\deropx f(x) = \frac{1}{(I - \lop)x}\left( I - \lop \right)f(x)
\end{equation}
which allows us to formally invert the operator $\deropx$ to get the
inverse operator $\derinv$ through the definition
\begin{equation}
\derinv = \sum_{j=0}^{\infty}{\lop}^j \left\{ (I - \lop)x \right\}
\end{equation}
As an illustrative example, let us compute the operation of $\derinv$ on
$f(x)=x^n$ (say). We then have
\beq
\derinv \  x^n
&=& \sum_{j=0}^{\infty}{\lop}^j \left\{ (I - \lop)x \right\}x^n\nonum \\
&=& \sum_{j=0}^{\infty}{\lop}^j \left\{ x^{n+1} - (\lop
x) x^n\right\}\nonum \\
&=& \sum_{j=0}^{\infty}
        \left\{
        {\lop}^j x^{n+1} - ({\lop}^{j+1}x)({\lop}^{j}x^n)
        \right\}\label{inverse}
\eeq
The operation of $\deropx$ on (\ref{inverse}) gives
\bseq
\deropx \cp \derinv x^n
&=& \deropx \sum_{j=0}^{\infty}
        \left\{
        {\lop}^j x^{n+1} - ({\lop}^{j+1}x)({\lop}^{j}x^n)
        \right\} \\
&=& \left(\frac{1}{(\lop x) - x}\right)
    \Bigg\{
    \lop \sum_{j=0}^{\infty}
        \left\{
        {\lop}^j x^{n+1} - ({\lop}^{j+1}x)({\lop}^{j}x^n)
        \right\} \\
& &  \mbox{\hspace{25mm}} -   \sum_{j=0}^{\infty}
        \left\{
        {\lop}^j x^{n+1} - ({\lop}^{j+1}x)({\lop}^{j}x^n)
        \right\}
    \Bigg\} \\
&=&  \left(\frac{1}{(\lop x) - x}\right)
    \Bigg\{
      \sum_{j=0}^{\infty}
        \big\{
        {\lop}^{j+1} x^{n+1} - ({\lop}^{j+2}x)({\lop}^{j+1} x^n)\\
& &  \mbox{\hspace{25mm}}      -
        {\lop}^j x^{n+1} + ({\lop}^{j+1}x)({\lop}^{j}x^n)
        \big\}
     \Bigg\} \\
&=&  \left(\frac{1}{(\lop x) - x}\right)
     \Bigg\{
        \lop x^{n+1} - ({\lop}^2 x)(\lop x^n)
        - x^{n+1} + (\lop x)x^n +\\
& &  \mbox{\hspace{25mm}}   \lop ^2 x^{n+1} -(\lop ^3 x)(\lop ^2 x^n)
        - (\lop x^{n+1}) + (\lop ^2 x)(\lop x^n)     + \cdots
      \Bigg\} \\
&=& \left( \frac{1}{( \lop x) - x} \right)
        \left\{ ( \lop x)x^n - x^{n+1} \right\} \\
&=& x^n
\eseq
Hence the necessary condition
\begin{equation}
\deropx \cp \! \derinv = I
\end{equation}
is satisfied by the inverse derivative $\derinv$.

\subsection{Eigenfunctions}

The eigenfunctions of the deformed derivative operator $\deropx$ can be
shown to admit a product representation. The following proposition
justifies this statement.
\bpr
If
\begin{equation}\label{elameig}
  \deropx \Elam(x) = \Elam(x)
\end{equation}
then the function $\Elam(x)$ is given by
\begin{equation}
\Elam (x) = \prod_{j=0}^{\infty}\frac{1}{1+{\lop}^j (\lop -I)x}
\end{equation}
\epr
{\bf Proof:}
Let the function $\Elam(x)$ be defined as
\begin{equation}\label{lamexp}
  \Elam(x) = \prod_{j=0}^{\infty}c_j(x)
\end{equation}
where $c_j(x)$ is a sequence of functions of $x$
\begin{equation}\label{cks}
  c_j(x) = {\lop}^j c_0(x) = c_0({\lop}^j x)
\end{equation}
Then
\begin{eqnarray*}
\deropx \Elam (x)
&=& \left\{ \frac{1}{(\lop -I)x} \right\} \left( \lop - I \right)
        \prod_{j=0}^{\infty}c_j(x)\nonum \\
&=& \left\{ \frac{1}{(\lop -I)x} \right\}
     \left\{ \prod_{j=0}^{\infty} \lop c_j(x)
             - \prod_{j=0}^{\infty}c_j(x) \right\}\\
&=& \left\{ \frac{1}{(\lop -I)x} \right\}
     \left\{ \prod_{j=0}^{\infty} c_{j+1}(x)
             - \prod_{j=0}^{\infty} c_j(x) \right\}\\
&=& \left\{ \frac{1}{(\lop -I)x} \right\}
     \left\{ \prod_{j=0}^{\infty} \frac{c_{j}(x)}{c_0(x)}
             - \prod_{j=0}^{\infty} c_j(x) \right\}\\
&=& \left\{ \frac{1}{(\lop - I)x} \right\}
        \left\{ c_{0}^{-1}(x) - 1\right\}\prod_{j=0}^{\infty}c_j(x)\\
&=& \left\{ \frac{1}{(\lop -I)x} \right\}
                \left\{ c_{0}^{-1}(x) - 1 \right\} \Elam (x)
\end{eqnarray*}
Since the function $c_0(x)$ (cf. (\ref{cks}))is an arbitrary function of
$x$, we may choose
\begin{equation}
  c_{0}^{-1}(x) = 1 + (\lop x) - x
\end{equation}
so that the functions $c_j(x)$ can be cast in the form
\bseq
c_j(x)
&=& {\lop}^j c_0(x) \\
&=& c_0({\lop}^j x) \\
&=& \frac{1}{1 + {\lop}^j (\lop -I)x}
\eseq
which, in turn, gives the desired eigenfunction
\begin{equation}\label{lexp}
  \Elam (x) = \prod_{j=0}^{\infty}\frac{1}{1+{\lop}^j (\lop -I)x}
\end{equation}
In view of the relation $\deropx \Elam(x) = \Elam(x)$, the eigenfunction 
$\Elam(x)$ can be interpreted as the $\Omega$--exponential function. An
alternative description of the eigenfunctions $\Elam(x)$ proceeds via the
series representation in terms of the inverse operator $\derinv$. For this
purpose, we write
\begin{equation}\label{series}
\Elam(x) = \sum_{j=0}^{\infty} e_j(x)
\end{equation}
where functions $e_j(x)$ are defined as
\begin{equation}
e_j(x) = \deropx e_{j+1}(x) \qquad \text{or} \qquad
\derinv e_j(x) = e_{j+1}(x)
\end{equation}
with the normalisation $e_0(x) = 1$. Relation (\ref{series}) now becomes
\begin{equation}
\Elam(x) = \sum_{j=0}^{\infty} ({\derinv})^j e_0(x)
=\sum_{j=0}^{\infty} ({\derinv})^j
\end{equation}

\section{Examples of the Operator $\lopr$}

As remarked in Sec. $2$, the linear operator $\lop$ is not only
multiplicative over functions of $x$, but is also characterised by a
deformation parameter $\lambda$. For a complete specification, however, it
carries an additional index $k\in \mathbb{Z}^+$ signifying the
degree of the operator $\lop$. We illustrate this by considering a class
of operators given by
\begin{equation}\label{def1}
\lopr = \exp\left({-\lambda x^{k+1}{\partial}_x}\right)
\end{equation} 
and study their action on functions of the real variable $x$.
\bpr
If $f(x)=x^n$, $(n\in \mathbb{Z}^+)$, then
\begin{equation}\label{prop2}
\lopr  x^n = \left\{
                \frac{x}{\left( 1+\lambda k x^k \right)^{\frac{1}{k}}}
            \right\}^n
\end{equation}
\epr
{\bf Proof:}
Let us consider the operation
\beq
\exp\left({-\lambda x^{k+1}{\partial}_x}\right)x^n
&=& \sum_{n=0}^{\infty} \frac{(-\lambda)^n}{n!} \left( x^{k+1}\dx
\right)^n x^n\nonum \\
&=&       x^n + (-\lambda) n x^{k+n}
        + \frac{{(-\lambda)}^2}{2!} n (n+k)x^{2k +n}
        + \cdots \nonum \\
&=&     x^n \left\{
          1 + \left(\frac{-n}{k}\right) (\lambda k x^k)
            + \frac{1}{2!}
            \left(\frac{-n}{k}\right)\left(\frac{-n}{k}-1\right)
            (\lambda k x^k)^2
            + \cdots
            \right\} \nonum \\
&=&  x^n \left( 1+ \lambda k x^k \right)^{-\frac{n}{k}}\nonum \\
&=&  \left\{ \frac{x}{\left( 1+ \lambda k x^k \right)^{\frac{1}{k}}}
     \right\}^n \label{proof2}
\eeq
On the other hand, if we consider another class of operators
\begin{equation}\label{def2}
\lopr = \exp\left({
                -\frac{\lambda}{2}\{\xop^{k+1},{\partial}_x}\}
            \right)
\end{equation}
where $\{\ , \ \}$ is the anticommutator, then their operation given by
the adjoint action on the position operator $\xop$ results in the
equivalence of (\ref{def1}) and (\ref{def2}) which can be readily seen as
follows:
\beq
\lopr \xop \loprinv
&=& \exp\left(  -\frac{\lambda}{2}\{ \xop^{k+1},{\partial}_x \}
        \right)
    \xop
    \exp\left(   \frac{\lambda}{2}\{ \xop^{k+1},{\partial}_x \}
        \right) \\
&=&     \xop
      + \left[
        - \frac{ \lambda}{2} \{ \xop^{k+1},{\partial}_x \}, \xop
        \right] \nonum \\
& &   \mbox{\hspace{5mm}}
      + \frac{1}{2!}
        \left[
        -\frac{\lambda}{2} \{ \xop^{k+1},{\partial}_x \},
        \left[
        -\frac{\lambda}{2} \{ \xop^{k+1},{\partial}_x \}, \xop
        \right]
        \right]
      + \cdots
\eeq
using the Baker-Campbell-Hausdorff formula. The evaluation of the nested
commutators proceeds through the following steps
\bseq
\left[
       - \frac{ \lambda}{2} \{ \xop^{k+1},{\partial}_x \}, \xop
\right]
&=&   - \frac{ \lambda}{2}
        \left[ \xop^{k+1}\dx + \dx \xop^{k+1} , \xop \right] \\
&=&   - \lambda \xop^{k+1} \\
&=&     \xop
        \left( -\frac{1}{k} \right)\left(\lambda k \xop^k
        \right)\\
& &     \\
\left[
      -\frac{\lambda}{2} \{ \xop^{k+1},{\partial}_x \},
        \left[
         -\frac{\lambda}{2} \{ \xop^{k+1},{\partial}_x \}, \xop
        \right]
\right]
&=&     \left[ -\frac{\lambda}{2} \{ \xop^{k+1},{\partial}_x \},
                - \lambda \xop^{k+1}
        \right] \\
&=&     {\lambda}^2 (k+1)\xop^{2k+1} \\
&=&     \xop \left( -\frac{1}{k} \right)\left(-\frac{1}{k}-1\right)
        \left(\lambda k \xop^k \right)^2
\eseq
and so on. All these terms add up to give
\beq
\lopr \xop \loprinv
&=& \xop
        \left\{ I
        + \left( -\frac{1}{k} \right)\left(\lambda k \xop^k \right)
        + \frac{1}{2!} \left( -\frac{1}{k}
\right)\left(-\frac{1}{k}-1\right)
          \left(\lambda k \xop^k \right)^2
        + \cdots
        \right\}\nonum \\
&=&  \frac{\xop}{\left( I + \lambda k \xop^k \right)^{\frac{1}{k}}}
\eeq
leading finally to the result
\begin{equation}
\lopr \xop^n \loprinv
= \left\{
        \frac{\xop}{\left( 1 + \lambda k \xop^k \right)^{\frac{1}{k}}}
  \right\}^n
\end{equation}
which is precisely the same as (\ref{prop2}). We, therefore, conclude that
the class of operators (\ref{def1}) and  (\ref{def2}) are essentially
equivalent. 
\par
A cursory inspection of the operators $\lopr$ ((\ref{def1}) or
(\ref{def2})) helps in defining the relations
\begin{equation}
\left( \lopr \right)^{-1} = {\hat{\Omega}}_{ -\lambda , k}
\end{equation}
and
\begin{equation}
\lopr \cp {\hat{\Omega}}_{ \mu , k} = {\hat{\Omega}}_{\lambda + \mu , k}
\end{equation}
These exhibit interesting consequences in that if $\lopr x=y$,
then $\loprinv y=x$. Writing explicitly, if we have
\begin{equation}
y = \lopr x = \frac{x}{\left( 1+ \lambda k x^k \right)^{\frac{1}{k}}}
\end{equation}
then
\begin{equation}
x = \loprinv y
  = \frac{y}{\left( 1 - \lambda k y^k \right)^{\frac{1}{k}}}
\end{equation}
where we have made use of (\ref{proof2}) for the case when $n=1$. This is
an interesting off-shoot of the foregoing considerations, allowing
construction of a class of other operators possessing this unique 
property.
\par
In more general terms, let us consider the operator
\begin{equation}
\Lmn = \exp\left( \lambda f(x)\dx \right)
\end{equation}
such that $\lopf x = y$. In view of the above structures, the functions
$f(x)$ and $f(y)$ are invertible and therefore admit the forms
\begin{equation}
f(y)  = \frac{f(x)}{1+\lambda f(x)}\, ,
\qquad \qquad
f(x)  = \frac{f(y)}{1-\lambda f(y)}
\end{equation}
For instance, if we compute the action of $\Lmn$ on $f(x) = \exp(-\mu x)$,
say, we find
\begin{equation}
\Lmn x = x + \mu^{-1} \ln \left[1 + \lambda \mu e^{-\mu x} \right]
\label{linvert}
\end{equation}
which constitutes a two--parameter $(\lambda,\mu)$ deformation of
the derivative operator. Note, however, that $\Lmn \ra I$ in the
limiting case when $\lambda \ra 0$ which is in complete conformity with
what we expect. In a similar way, we also evaluate the limit of $\Lmn x$
(\ref{linvert}) when $\mu \ra 0$ with the result
\begin{equation} 
\mbox{lim}_{\mu \ra 0}
        \left\{
           x + \mu^{-1} \left(\lambda \mu e^{-\mu x} \right)
             + O(\mu)
        \right\} = x + \lambda
\end{equation}
This is obviously a manifestation of the finite difference operator
$\mbox{lim}_{\mu \ra 0} \Lmn = \exp(\lambda \dx)$, characterised by the
parameter $\lambda$. We remark in passing that the form 
(\ref{linvert}) can be alternatively represented as
\begin{equation}
\Lmn x = \mu^{-1} \ln \left[e^{\mu x} +\lambda \mu \right]
\end{equation}

\subsection{The associated derivative}

The operation of the deformed derivative associated with $\lopr$ on a
function $f(x)$ can be expressed as
\begin{equation}
\derlr f(x) = \frac{f(\lopr x) - f(x)}{(\lopr x) - x}
\end{equation}
Now
\beq
(\lopr x) - x
&=&  \frac{x}{\left( 1+ \lambda k x^k \right)^{\frac{1}{k}}} -x
     \nonum   \\
&=&  x \left\{
        \frac{ 1 - (1+\lambda k x^k)^{\frac{1}{k}}}
                  {(1+\lambda k x^k)^{\frac{1}{k}}}
       \right\}
\eeq
so that
\begin{equation}
\derlr f(x)
=       \frac {(1+\lambda k x^k)^{\frac{1}{k}}}
		{ x\left[ 1 - (1+\lambda k x^k)^{\frac{1}{k}}\right] }
       \left\{
       f\left(
              \frac{x}{\left( 1+ \lambda k x^k \right)^{\frac{1}{k}}}
        \right)
       - f(x)
       \right\}
\end{equation}
If we choose $f(x) = x^n$, then
\begin{equation}
\derlr x^n =  [[n]]_{\lambda , k}\, x^{n-1}
\end{equation}
where
\beq
[[n]]_{\lambda , k} \, = \, [[n]]_{\lambda , k}(x)
&=& \frac{
        \left\{
        \frac{1}{\left( 1+ \lambda k x^k \right)^{\frac{1}{k}}}
        \right\}^n -1}
        {\left\{
        \frac{1}{\left( 1+ \lambda k x^k \right)^{\frac{1}{k}}}
        \right\} -1} \\
&=& \sum_{n=0}^{\infty} \left( 1+ \lambda k x^k \right)^{-\frac{n}{k}}
\eeq
Note that $[[n]]_{\lambda , k}$ is now a function of $x$.

\subsection{Deformed inverse eigenfunctions}

We have already seen that the $\Omega$--exponential functions $\Ellam(x)$
are the eigenfunctions of $\deropx$. This enables us to write down the
eigenvalue equation
\begin{equation}\label{evalue}
\deropx \Ellam(\mu x) = \mu \Ellam(\mu x) 
\end{equation}
where $\mu$ is regarded as an eigenvalue of $\deropx$ with the
corresponding eigenfunction $\Ellam(\mu x)$. Explicitly, this can be
written as
\begin{equation}\label{Emu}
  \lopr \Ellam(\mu x) = \left\{1 + \mu (\lopr -I)x \right\} \Ellam(\mu x)
\end{equation}
which reduces to 
\begin{equation}\label{ev2}
\lopr \Ellam( x) = \left\{1 +  (\lopr -1)x \right\} \Ellam( x)
\end{equation}
for $\mu =1$. Rewriting (\ref{Emu}), we have
\beq
\Ellam(\mu x)
&=& \loprinv \left\{ 1 + \mu (\lopr -1)x \right\} \Ellam(\mu x) \\
&=& \left\{ 1 + \mu (1 - \loprinv )x \right\}
        \Ellam (\mu \loprinv x)\\
&=&	\left\{ 1 + \mu (1 - \hat{\Omega}_{-\lambda, k} )x \right\}
        \Ellam (\mu \hat{\Omega}_{-\lambda,k}\, x) \label{ev1}
\eeq
where use has been made of the relation $\loprinv = \hat{\Omega}_{-\lambda
, k}$. Letting $\lambda \ra -\lambda$ and $\mu = -1$ in (\ref{ev1}),
we obtain
\beq
E_{-\lambda}(-x)
&=&     \left\{ 1 + ( \hat{\Omega}_{\lambda , k} - 1 )x \right\}
        E_{-\lambda}(-\hat{\Omega}_{\lambda , k}\, x) \\
&=&     \left\{ 1 + ( \hat{\Omega}_{\lambda , k} - 1 )x \right\}
        \hat{\Omega}_{\lambda,k}\,E_{-\lambda} (-x) \label{-Emu}
\eeq
Equations (\ref{ev2}) and (\ref{-Emu}) lead to
\begin{equation}
\lopr \Ellam (x)\hat{\Omega}_{\lambda , k}\,E_{-\lambda}(-x) =
\Ellam (x)E_{-\lambda}(-x)
\end{equation}
i.e.
\begin{equation}
\lopr \left\{ \Ellam (x) E_{-\lambda}(-x) \right\} = \Ellam
(x) E_{-\lambda}(-x)
\end{equation}
or
\begin{equation}
(\lopr - I) \left\{ \Ellam (x) E_{-\lambda}(-x) \right\} = 0
\end{equation}
Hence
\begin{equation}
\derlr \left\{ \Ellam (x) E_{-\lambda}(-x) \right\} = 0
\end{equation}
which clearly shows that the product $\Ellam (x) E_{-\lambda}(-x)$
under the action of the deformed derivative operator $\derlr$
differentiates to zero. Ordinarily, we would expect that only
constant functions are endowed with this property. An important
consequence of this observation is that any function $F(x)$, not
necessarily a constant, with the property $F(x)=F(\lop x)$ will
differentiate to zero under the action of the generalised derivative
$\derlr$.
\par
We now consider the case where $\lopr$ is identified with the translation
operator, i.e., $\lopr = \hop = e^{h \dx}$. It is then easy 
to see from (\ref{ev2}) that
\beq
  \hop E_h(x)
  &=& \left\{ 1 + (\hop - 1)x \right\} E_h(x)\nonum \\
  &=& (1+h)E_h(x)
\eeq
Consequently, the eigenfunctions of the translation operator $\hop$ are
exactly the same as those of the deformed derivative operator associated
with $\hop$. On the other hand, if $\lopr$ is taken to be the
$q$--dilation operator, i.e., $\lopr = \qop = \exp(\ln (q) x \dx )$, then
we obtain
\begin{equation}
E_q^{-1}(x) = E_{q^{-1}}(-x) 
\end{equation}
as expected.

\section{The Conformal (M\"obius) Transformation}

The mathematical framework developed thus far will now be cast in a proper
perspective so as to make contact with situations that hold promise for
possible applications. The translation operator $\hop$ and the
$q$--dilation operator $\qop$, which give rise to the finite difference
derivative and the $q$--derivative, respectively, may be regarded as
generators of a one-dimensional affine group. In view of our consideration
of the deformed derivative operator, of which $\hop$ and $\qop$ are
special cases, it would seem instructive to visualise a situation where
the operator $\lopr$, affecting the generalisation, can indeed be assigned
the role of an additional generator, resulting in some larger algebraic
structure. To this end, we take recourse to the special conformal
(M\"obius) transformation and carry out an embedding in conjunction with
the generators $\hop$ and $\qop$.
\par
The action of a matrix $M = \left( \begin{smallmatrix}a & b\\c &
d\end{smallmatrix} \right)$ on a variable $x$ is defined as
\begin{equation}
M \triangleright x = \frac{ax + b}{cx + d}, \qquad
ad\neq bc
\label{affine}
\end{equation}
The operator $\lopr$ (cf. (\ref{def1})) for $k=1$ takes the form
\begin{equation}
\lop = \hat{\Omega}_{\lambda , k=1} = \exp(-\lambda x^2\dx)
\end{equation}
whose action on a variable $x$ (cf. (\ref{prop2})) is given by
\begin{equation}
\lop x = \frac{x}{1+\lambda x}
\end{equation}
Defining the $\lam$, $Q$ and $H$ matrices as
\begin{equation}
\lam = \left( \ba{cc} 1&0 \\ \lambda &1 \ea \right)\, , \qquad
   Q = \left( \ba{cc} q^{\frac12}&0 \\ 0& q^{-\frac12} \ea \right)\, ,
\qquad H = \left( \ba{cc} 1&h \\ 0&1 \ea \right)
\end{equation} 
the corresponding actions on $x$ are represented as
\begin{equation} 
\lam
\triangleright x = \frac{x}{\lambda x + 1}\, , \qquad
   Q \triangleright x = q x\, , \qquad
   H \triangleright x = x+h
\end{equation}
The general M\"obius operator $M$ is then obtained via the composition
\beq
M &=& \lam \, Q \, H \\
\text{i.e.}\qquad \left( \ba{cc} a&b \\ c&d \ea \right)
&=& \left( \ba{cc} q^{\frac12}& q^{\frac12}h \\
           q^{\frac12}\lambda & q^{\frac12}\lambda h + q^{-\frac12}
           \ea
     \right)
\eeq
Equivalently, we express the M\"obius transformations on the variable $x$
in the form 
\beq
\hat{M}\, x &=& \lop \, \qop \, \hop \, x \\ 
	&=& \exp \left[-\lambda x^2\dx + \ln(q)x\dx + h\dx\right] \, x
\eeq
It turns out that $ad-bc=1$. Thus the set of $2\times 2$ matrices $M$
generate the $SL(2)$ group structure. In other words, the operator $\lop$
serves as an extra generator, in addition to the other two generators
typified by the translation operator $\hop$ and the $q$--dilation operator
$\qop$, to generate the group $SL(2)$ as well as its subgroups. The
operator $\lop$ is designated as the reduced M\"obius transformation
operator and conforms to the realisation of the differential operators
in view of the transformation (\ref{affine}).

\section{Concluding Remarks}
A modest attempt has been made to construct a generalised version of the
deformed derivative operators within the context of a linear operator
$\lop$ whose action on functions of a real variable is considered to
be multiplicative. The veracity of the mathematical formulation is checked
against suitable illustrative examples. We note that the familiar forms of
the finite difference and the $q$--derivatives are readily obtained as
particular cases of the generalised structure. Such an observation
lends credence to the validity of the procedure adopted in carrying
out this construction. We have also addressed, though briefly, the
possibility of casting the operator $\lop$ in a format compatible with the
special conformal or reduced M\"obius transformation operator. The 
inclusion of this operator, along with the translation operator $\hop$ and
the $q$--dilation operator $\qop$, to generate the $SL(2)$ group is
expected to be significant in seeking further applications to situations
of physical interest such as the conformal transformations. The
construction also forms a natural paradigm to investigate generalisation
of the $q$--oscillators and the $q$--calculus.

\section*{Acknowledgements}
The first author is grateful to Prof. E. Zeidler for the kind hospitality
at the Max Planck Institute for Mathematics in the Sciences, Leipzig,
where this work was carried out. The second author is supported by a
research fellowship of the Royal Commission for the Exhibition of 1851.

\end{document}